\documentclass[a4, 10pt]{amsart}
\usepackage{amssymb}
\usepackage{amstext}
\usepackage{amsmath}
\usepackage{amscd}
\usepackage{latexsym}
\usepackage{amsfonts}
\usepackage[dvips]{graphicx}

\newtheorem{thm}{Theorem}[subsection]
\newtheorem{proposition}[thm]{Proposition}
\newtheorem{lem}[thm]{Lemma}
\newtheorem{corollary}[thm]{Corollary}

\theoremstyle{definition}
\newtheorem{dfn}{Definition}
\newtheorem{example}[thm]{Example}

\theoremstyle{remark}
\newtheorem{rmk}{Remark}

\numberwithin{equation}{section}

\begin{document}

\title[On the Jordan decomposition of tensored matrices]{On the Jordan decomposition of tensored matrices of Jordan canonical forms}
\author{Kei-ichiro IIMA}
\address{Graduate School of Natural Science and Technology, Okayama University, 1-1, Naka 3-chome, Tsushima, Okayama 700-8530, JAPAN}
\email{iima@math.okayama-u.ac.jp}
\author{Ryo IWAMATSU}
\address{Graduate School of Education, Saitama University, Shimo-Okubo 255, Sakura, Saitama, Saitama 338-8570, JAPAN}
\email{05AC201@post.saitama-u.ac.jp}
\thanks{
\endgraf
}

\maketitle
\begin{abstract}
Let $k$ be an algebraically closed field of characteristic $p \ge 0$.
We shall consider the problem of finding out a Jordan canonical form of \ \ \ $J(\alpha,\,s) \otimes_{k} J(\beta,\,t)$, where $J(\alpha,\,s)$ means the Jordan block with eigenvalue $\alpha \in k$ and size $s$.
\end{abstract}
\section{Introduction} 
To construct graded local Frobenius algebras over an algebraically closed field $k$,
it is important to find out a Jordan canonical form (simply, JCF) of tensor product of square matrices.
In fact, it is known that 
any graded local Frobenius algebra is of the form of $\Lambda(\varphi,\,\gamma)=T(V)/R(\varphi,\,\gamma)$,
where $V$ is a finite dimensional $k$-vector space,
$\gamma$ an element of $GL(V)$,
and $\varphi : V^{\otimes n} \to k$ a $k$-linear map satisfying several conditions.
Further, if we decompose as $(V,\,\gamma)=\bigoplus_{i}(V_{i},\,\gamma_{i})$,
then the conditions of $\varphi$ can be described in terms of
each $\varphi_{i_{1}\dots i_{r}} : V_{i_{1}} \otimes \cdots \otimes V_{i_{r}} \to k$.
Then, we have to consider a JCF of $\gamma_{i_{1}} \otimes \cdots \otimes \gamma_{i_{r}}$
as an element in $GL(V_{i_{1}} \otimes \cdots \otimes V_{i_{r}})$.
(For detail, refer to T. Wakamatsu \cite{Wak}).

Let $k$ be an algebraically closed field of characteristic $p \ge 0$, and $J(\alpha,\, s), J(\beta,\,t)$ Jordan blocks over $k$.
We shall consider the problem of finding out a JCF of $J(\alpha ,s) \otimes J(\beta ,t)$, where $\otimes$ means $\otimes_{k}$ $(s \le t)$.

Over an algebraically closed base field of characteristic zero, this problem has been solved by many authors including T. Harima and J. Watanabe \cite{HW}, and A. Martsinkovsky and A. Vlassov \cite{MV} etc.
M. Herschend \cite{H} solve it for extended Dynkin quivers of type $\tilde{\Bbb A}\sb n$, with arbitrary orientation and any $n$.
In this note we solve it for any characteristic $p \geq 0$. That is, we obtain two ways to determine the Jordan decomposition of the tensored matrix $J(\alpha, s) \otimes J(\beta, t)$.

In the case of $\alpha \beta = 0$, the tensored matrix $J(\alpha ,s) \otimes J(\beta ,t)$ has the same direct sum decomposition as in Theorem \ref{MV} independently of characteristic of the base field $k$ in Proposition \ref{ab=0}.
In the case of $\alpha \beta \neq 0$, our problem is reduced to the problem of finding the indecomposable decomposition of $R$ as a $k[\theta]$-module, where $R$ means the quotient ring $k[x,y]/(x^{s},y^{t})$, $\theta = x+y$ and $k[x,y]$ be a polynomial ring over $k$ . 
In the section 2.1, we regard finding the indecomposable decomposition of $R$ as calculating the partition $\bold{c} = (c_{1}, c_{2}, \ldots , c_{r})$ of $st$ in Lemma \ref{lemma}.
Then, we are able to determine the Jordan decomposition of tensored matrix $J(\alpha ,s) \otimes J(\beta ,t)$.
In the section 2.2, we show another algorithm.
The idea is finding out elements that determine the indecomposable decomposition of $R$ as a $k[\theta]$-module.
In Theorem \ref{3}, we show that we can find out $s$ homogeneous elements $\omega_{0}, \omega_{1}, \dots, \omega_{s-1}$ of $R$ such that $R = \bigoplus_{i=0}^{s-1} k[\theta] \omega_{i}$, where the degree of $\omega_{i}$ is $i$ for each $0 \le i \le s-1$.
And applying this result, we show an algorithm for computing a JCF of $J(\alpha,\,s) \otimes J(\beta,\,t)$ in Theorem \ref{15}.




\section{Main results}

Throughout this section, let $k$ be an algebraically closed field.
For an integer $s\ge 1$ and an element $\alpha\in k$, let
$$
J(\alpha ,s)=
\begin{pmatrix}
\alpha & 1 \\
 & \ddots & \ddots \\
 & & \alpha & 1 \\
 & & & \alpha
\end{pmatrix}
$$
denote the Jordan block of size $s\times s$ with an eigenvalue $\alpha$.

\begin{thm}\cite[Theorem 2]{MV}\label{MV}
Suppose that $k$ has characteristic zero.
Then the following holds for integers $s\le t$ and $\alpha, \beta\in k$:
$$
J(\alpha ,s)\otimes J(\beta ,t)=
\begin{cases}
J(0 ,s)^{\oplus t-s+1}\oplus\bigoplus _{i=1}^{2s-2}J(0,s-\lceil\frac{i}{2}\rceil) & \text{if}\quad\alpha =0=\beta \\
J(0 ,s)^{\oplus t} & \text{if}\quad\alpha = 0 \neq \beta \\
J(0 ,t)^{\oplus s} & \text{if}\quad\alpha \neq 0 = \beta \\
\bigoplus _{i=1}^{s} J(\alpha\beta ,s+t+1-2i) & \text{if} \quad \alpha \neq 0 \neq \beta
\end{cases}
.
$$
\end{thm}

\begin{rmk}
If one of the eigenvalues $\alpha$ and $\beta$ equals zero, then the tensored matrix $J(\alpha ,s) \otimes J(\beta ,t)$ has the same direct sum decomposition as in Theorem \ref{MV} independently of characteristic of the base field $k$ (Proposition \ref{ab=0}).
\end{rmk}

\begin{thm}\label{main}
There is an algorithm to determine the Jordan decomposition of the tensored matrix $J(\alpha ,s) \otimes J(\beta ,t)$, which has  an independent description of the characteristic of the base field $k$.
\end{thm}

\begin{rmk}

(1)The matrix $J(\alpha ,s)$ represents the action of $X$ on $k[X]/(X - \alpha)^{s}$ as a $k[X]$-module.

(2)The tensored matrix $J(\alpha ,s) \otimes J(\beta ,t)$ is triangular. Therefore its eigenvalue is $\alpha \beta$.

(3)One has an isomorphism $$k[X]/(X - \alpha)^{s} \otimes k[Y]/(Y - \beta)^{t} \cong k[X,Y]/((X - \alpha)^{s},(Y - \beta)^{t})$$ of $k$-algebras.
\end{rmk}

Tensored matrix $J(\alpha,s)\otimes J(\beta,t)$ represents the action of $XY$ on $k[X,Y]/((X-\alpha )^s,(Y-\beta )^t)$ as a $k[XY]$-module.

\subsection{The method for calculating numerical values of tensored matrices.}

\begin{lem}\label{lemma}
Put $R=k[X,Y]/((X-\alpha )^s,(Y-\beta )^t)$, which we regard as a $k[Z]$-module through the map $k[Z]\to R$ given by $Z \mapsto XY$.
Then there is a sequence of integers such that $c_1\ge c_2\ge \cdots\ge c_r\ge 1$
$$
R\cong \bigoplus _{i=1}^r k[Z]/(Z-\alpha\beta )^{c_i}
$$
of $k[Z]$-modules.
\end{lem}

This means that $J(\alpha ,s) \otimes J(\beta ,t) = \bigoplus _{i=1}^{r} J(\alpha \beta ,c_{i})$. 
We can regard $\bold{c}=(c_{1},c_{2}, \ldots ,c_{r})$ as a partition of $st$ in obvious manner. 
The main problem is to determine the partition $\bold{c}$.
For this purpose let $\bold{b}=(b_{1},b_{2}, \ldots ,b_{r'})$ be the partition conjugate to $\bold{c}$. 
Put $z=Z- \alpha \beta$. 
Note that $b_{i}= \# \{ j | c_{j} \ge i \} = \dim _{k}(z^{i-1}R/z^{i}R)$. 
Setting $a_{i}= \dim _{k}(R/z^{i}R)$, we have $b_{i}=a_{i}-a_{i-1}$. 
Therefore, it is sufficient that we calculate the value of $a_{i}$ for each case.

If one of the eigenvalues $\alpha$ and $\beta$ equals zero, then the result is independent of the characteristic of $k$ as we show in the next proposition.
\begin{proposition}\label{ab=0}
We have the following equalities;
$$
a_{i}=
\begin{cases}
(s+t)i-i^{2} \ \ (1 \le i \le s) & \text{if}\quad\alpha =0=\beta \\
ti \ \ (1 \le i \le s) & \text{if}\quad\alpha = 0 \neq \beta \\
si \ \ (1 \le i \le t) & \text{if}\quad\alpha \neq 0 = \beta
\end{cases}
.
$$

Therefore we get 

$$
J(\alpha ,s)\otimes J(\beta ,t)=
\begin{cases}
J(0 ,s)^{\oplus t-s+1}\oplus\bigoplus _{i=1}^{2s-2}J(0,s-\lceil\frac{i}{2}\rceil) & \text{if}\quad\alpha =0=\beta \\
J(0 ,s)^{\oplus t} & \text{if}\quad\alpha = 0 \neq \beta \\
J(0 ,t)^{\oplus s} & \text{if}\quad\alpha \neq 0 = \beta
\end{cases}
.
$$
\end{proposition}

\begin{proof}
Put $x=X-\alpha$ and $y=Y-\beta$.

(1) The case $\alpha =0= \beta$:\\
Since $R/z^iR=k[x,y]/(x^s,y^t,(xy)^i)$, we have $a_i= (s+t)i-i^{2}$.

(2) The case $\alpha =0 \neq \beta$:\\
Since $R/z^iR=k[x,y]/(x^s,y^t,x^i)$, we have $a_i=ti$.

(3) The case $\alpha \neq 0 = \beta$:\\
Since $R/z^iR=k[x,y]/(x^s,y^t,y^i)$, we have $a_i=si$.
\end{proof}

In the case of $\alpha \neq 0 \neq \beta$, then we have the following isomorphism of $k$-algebras, 
given by $X \mapsto x + \alpha$, $Y \mapsto \frac{- \alpha \beta}{y - \alpha}$:
$$
k[X,Y]/((X- \alpha)^{s},(Y- \beta)^{t},(XY- \alpha \beta)^{u}) \cong k[x,y]/(x^{s},y^{t},(x+y)^{u}).
$$
Using this isomorphism together with \cite[Proposition 4.4]{HMNW}\cite[Proposition 8]{HW}, we have the following proposition in the case of characteristic zero.

\begin{proposition}
Suppose that $\alpha \neq 0 \neq \beta$ and that $k$ has characteristic zero. 
Then we have
$$
\bold{b} =(\underbrace{s,s,\ldots ,s}_{t-s+1},s-1,s-1,s-2,s-2,\ldots ,1,1).
$$
Therefore we get $J(\alpha ,s)\otimes J(\beta ,t)=\bigoplus_{i=1}^{s} J(\alpha \beta ,s+t+1-2i)$.
\end{proposition}

\begin{proof}
Since the linear element $x+y \in k[x,y]/(x^s,y^t)$ is a strong Lefschetz element \cite{HW}. 
Namely, the multiplication map $\times (x+y)^{u} : k[x,y]/(x^s,y^t) _{i} \to k[x,y]/(x^s,y^t) _{i+u}$ is either injective or surjective,  for each $0 \le i \le s+t-2$. 
Then, we can easily compute $\dim _{k}(k[x,y]/(x^s,y^t,(x+y)^u))$ for each $1 \le u \le s+t-1$.
The assertion follows from this. 
\end{proof}

We consider in the rest the case where $\alpha \neq 0 \neq \beta$ and that $k$ is of positive characteristic $p$. 
Put $S=k[x,y]$, $R=k[x,y]/(x^{s},y^{t})$ and $A^{(u)}=R/(x+y)^{u}R$. 
To determine $a_{u} = \dim _{k}(A^{(u)})$, we may assume that $s \le t \le u$ without loss of generality. 
For each integer $u$ satisfying $s \le t \le u \le s+t-1$, we describe 
{\tiny
$$
(x + y)^{u} \equiv \binom{u}{s-1} x^{s-1} y^{u-s+1} + \binom{u}{s-2} x^{s-2} y^{u-s+2} + \cdots + \binom{u}{u-t+1} x^{u-t+1} y^{t-1} \pmod {(x^{s},y^{t})}.
$$
}
We set $q_{1}= \binom{u}{s-1}$, $q_{2}= \binom{u}{s-2}$,$\cdots$,
$q_{r}= \binom{u}{u-t+1}$ and $r = s + t - 1 - u$. 

We obtain the representation matrix of $R \stackrel{(x+y)^{u}}{\longrightarrow} R$ with respect to the natural base $\{ 1,x,y,x^{2},xy,y^{2},\ldots ,x^{s-1}y^{t-1} \}$ as follows;
$$
\begin{pmatrix}
 & & & & & & & \\
 & & & & & & & \\
H_{0} & & & & & & & \\
 & H_{1} & & & & & & \\
 & & H_{2} & & & & & \\
 & & & \ddots & & & & \\
 & & & & H_{r-2}& & & & \\
 & & & & & H_{r-1}& & & \\
\end{pmatrix},
$$
where
\begin{equation*}
H_{i} = 
\begin{pmatrix}
q_{i+1} & q_{i} & \cdots & q_{1}\\
q_{i+2} & q_{i+1} & \cdots & q_{2}\\
\vdots & \vdots & \ddots & \vdots \\
q_{r} & q_{r-1} & \cdots & q_{r-i}\\
\end{pmatrix}.
\end{equation*}

For each $0 \le i \le r-1$ the matrix $H_{i}$ is an $(r-i) \times (i + 1)$ matrix whose entries are integers.
We denote by $I_{i+1}(H_{i})$ the ideal of $\Bbb{Z}$ generated by $(i+1)$-minors of $H_{i}$ for $0 \le i \le r-1$. 
Obviously there exists an integer $\delta _{i} \ge 0$ such that $I_{i+1}(H_{i})$ = $\delta _{i} \Bbb{Z}$.
From the argument in the case of characteristic zero in \cite[Proposition 4.4]{HMNW}, we have $I_{i+1}(H_{i}) \otimes _{\Bbb{Z}} \Bbb{Q} \neq 0$, particularly $\delta _{i} \neq 0$, for any $0 \le i \le \lfloor (r-1)/2 \rfloor$.

\begin{proposition}
Under the same notation as above, for each $u$ satisfying $1 \le s \le t \le u \le s+t-1$, and for each $i$ satisfying $0 \le i \le \lfloor (r-1)/2 \rfloor (r=s+t-1-u)$, the following equalities hold;
{\small
\begin{equation*}
\delta_{i} = \gcd \{ S_{\lambda ^{j}}( \underbrace{1,1, \ldots ,1}_{u})| j=(j_{1},j_{2}, \ldots ,j_{i+1}) , 
1 \le j_{1} < j_{2} < \ldots < j_{i+1} \le r-i \},
\end{equation*}}
\noindent
where $\lambda ^{j}$ is the partition conjugate to $\mu ^{j} = (s-j_{1},s-j_{2}-1, \ldots ,s-j_{i+1}-i)$, and $S_{\lambda ^{j}}$ is the Schur polynomial.
\end{proposition}

\begin{proof}
Computation using Jacobi-Trudi formula \cite{FH} ,\cite{M}.
\end{proof}

Let
$$
0 \to S(-a) \oplus S(-b) \to S(-s) \oplus S(-t) \oplus S(-u) \stackrel{( x^{s} , y^{t} , (x+y)^{u} )}{\longrightarrow} S \to A^{(u)} \to 0
$$
be a minimal graded $S$-free resolution of $A^{(u)}$, where $1 \le s \le t \le u \le a \le b$. 
The Hilbert-Burch theorem implies that $a+b=s+t+u$, and the Hilbert series of $A^{(u)}$ is given as
$$
H_{A^{(u)}}(w)=
\frac {1 - w^{s} - w^{t} - w^{u} + w^{a} + w^{b}}{(1-w)^2}.
$$
It follows from this that $\dim _k(A^{(u)}) = st + su + tu - ab$. 
Letting $i_{0} = \min \{ i | \delta_{i} \equiv 0 \pmod{p} \}$, we get $a = u + i_{0}$ and $b = s + t - i_{0}$, since $a$ is the least value of degrees of relations of $(x^{s},y^{t},(x+y)^{u})$.
Thus, we can calculate the dimension of the $k$-vector space $A^{(u)}$, and hence the indecomposable decomposition of $J(\alpha ,s) \otimes J(\beta ,t)$.

\begin{thm}\label{alg}
We are able to compute a Jordan canonical form of $J(\alpha,\,s) \otimes J(\beta,\,t)$ by taking the following steps:
\begin{enumerate}
\item Every $\delta_{\bullet}$ is determined.
\item For each $1 \le u \le s+t-1$, $a_{u}$ is determined.
\item The partition $\bold{b}$ is determined.
\item The partition $\bold{c}$ is determined.
\item The Jordan decomposition of tensored matrix $J(\alpha ,s)\otimes J(\beta ,t)$ is determined.
\end{enumerate}
\end{thm}

From the discussion in Theorem \ref{alg}, one immediately obtains the following.

\begin{thm}
The tensored matrix $J(\alpha ,s) \otimes J(\beta ,t)$ has the same direct sum decomposition as in Theorem \ref{MV} if either $\mathrm{char}(k)\ge s+t-1$ or $I_{i+1}(H_{i}) \otimes _{\Bbb{Z}} k \neq 0$ for any $0 \le i \le \lfloor \frac{r-1}{2} \rfloor$.
\end{thm}


\subsection{The method for finding out elements that determine the indecomposable decomposition.}

$ $

In this subsection,
we show another algorithm for computing a JCF of $J(\alpha,\,s) \otimes 
J(\beta,\,t)$ 
via finding the indecomposable decomposition.
We have  already got the answer of our problem for case of $\alpha \beta =0$ 
by Proposition \ref{ab=0},
so we discuss only for case of $\alpha \beta \neq 0$.

We consider the indecomposable decomposition of $k[X,Y]/((X-\alpha)^{s},\, (Y-\beta)^{t})$ as a $k[XY]$-module.
As we stated in 2.1, we have an isomorphism $k[X,Y]/((X-\alpha)^{s},\, (Y-\beta)^{t}) \cong k[x,y]/(x^s,y^t)$.
Put $R=k[x,y]/(x^s,y^t)$, and $\theta=x+y$.
Thus, our problem is reduced to that of finding the indecomposable decomposition of $R$ as a $k[\theta]$-module.

It is clear that $R$ is a finite dimensional graded Artinian $k$-algebra.
So we write $R= \bigoplus_{i=0}^{s+t-2}R_{i}$.
And we immediately know that $\dim_{k} (R_{i})$ are written as $(1,\,2,\ldots ,s-1, \underbrace{s,\ldots ,s}_{t-s+1},s-1,\ldots ,1)$ for $0 \le i \le s+t-2$. We often use a figure for $R$ (Figure 1).

\begin{figure}
\centering
\includegraphics[clip, scale=0.5]{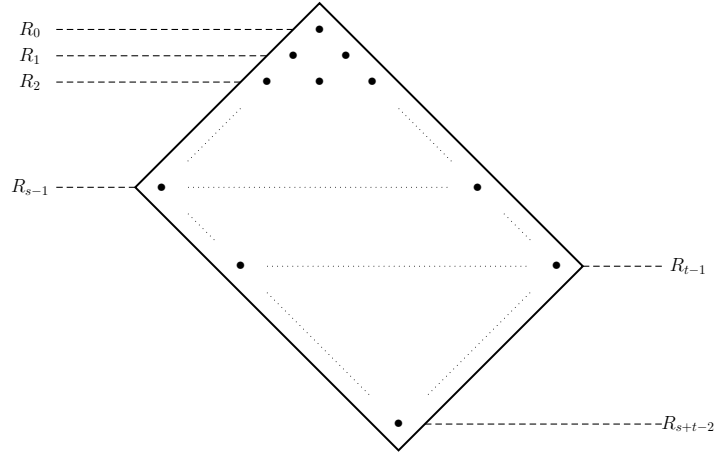}
\caption{An illustration of a basis of $R$. 
A bullet $\bullet$ stands for a base of $R$.}
\end{figure}

The subalgebra $k[\theta]$ of $R$ is \textit{uniserial}, which means 
that $k[\theta]$ has the only composition series as a $k[\theta]$-module.
We denote by $n$ the \textit{nilpotency} of $\theta$ (i.e. $\theta^{n} \neq 0$ and $\theta^{n+1}=0$).
And then, we can choose $\langle 1,\, \theta,\,\cdots,\, \theta^{n} \rangle$ as a $k$-basis of $k[\theta]$.
By easy calculation, we have the following inequality on $n$:

\begin{lem}\label{2}
We have $t-1 \le n \le s+t-2$. 
In particular, $n=s+t-2$ if $p=0$. 
\end{lem}

We describe the subalgebra $k[\theta]$ of $R$ in the figure for $R$ by drawing a polygonal line (Figure 2).

\begin{figure}
\centering
\includegraphics[clip, scale=0.5]{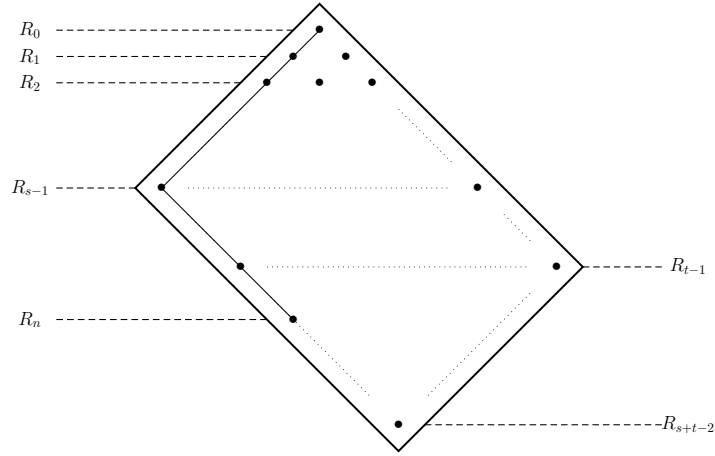}
\caption{An illustration of the subalgebra $k[\theta]$ of $R$. 
We consider the bullets on the polygonal line as $\langle 1,\, \theta,\,\cdots,\, \theta^{n} \rangle$.}
\end{figure}
 
Since the algebra 
$k[\theta]$ is uniserial,
any indecomposable summand $M$ of 
$R$ as a $k[\theta]$-module can be written as
$k[\theta] \omega$ for some element $\omega$ in $R$.
Hence we can write the indecomposable decomposition of 
$R$ as a $k[\theta]$-module such as:
\begin{equation}
R = \bigoplus_{i=1}^{r} k[\theta]\omega_{i} \>\>\>\> (\omega_{i} \in R). \tag{$*$}
\end{equation}

We shall call each element $\omega_{i}$ a \textit{generator} 
(for an indecomposable summand of $R$),
and the set $\{\omega_{1},\, \dots,\, \omega_{r}\}$, which consists of the generators in ($*$), 
a \textit{generating set} (for the indecomposable decomposition of $R$). 
A generating set is not unique.
However, we can prove the number of generators and 
that there exists the generating set which consists of homogeneous elements.

\begin{thm}\label{3}
There exists a generating set $\{\omega_{0},\, \omega_{1},\,\dots,\, \omega_{s-1}\}$
whose generator $\omega_{i}$ is an $i$-th degree homogeneous element. 
Hence,
$$
R = \bigoplus_{i=0}^{s-1} k[\theta]\omega_{i} \>\>\>\> (\omega_{i} \in R_{i}).
$$
\end{thm} 

In order to prove this theorem,
we have to prepare some lemmas and notations.

For a uniserial 
$k[\theta]$-submodule $M$ of $R$ generated by some homogeneous elements of $R$,
we denote by $\sigma(M)$ the \textit{socle degree} of $M$ as a $k[\theta]$-module.
In short, $\sigma(M)=d$ if $\mathrm{soc}_{k[\theta]}(M) \subseteq R_{d}$.
For example, $\sigma(k[\theta])=n$. And if $\theta^{n}x \neq 0$, then $\sigma(k[\theta]x)=n+1.$
The following lemma is checked easily:

\begin{lem}\label{4}
Let $\zeta,\,\eta$ be homogeneous elements of $R$.
If $\sigma(k[\theta]\zeta) \neq \sigma(k[\theta]\eta)$, 
then $k[\theta]\zeta \cap k[\theta]\eta =\{0\}$ holds. 
Hence $k[\theta]\zeta + k[\theta]\eta = k[\theta]\zeta \oplus k[\theta]\eta$.
\end{lem}

\begin{lem}\label{5}
Let $\kappa$ be a homogeneous element of $R$, and
put $d = \sigma(k[\theta]\kappa)$.
If  $t-1 \le d < s+t-2$,
then,
$$
\sum_{i=0}^{s+t-2-d} k[\theta]\kappa x^{i} 
= \bigoplus_{i=0}^{s+t-2-d} k[\theta]\kappa x^{i}.
$$
\end{lem}

\begin{proof}
Put $d^{\prime}=s+t-2-d$. And let the degree of $\kappa$ be $m$.
We now check $\theta^{d-m} \kappa x^{d^{\prime}} \neq 0$.
Since $\theta^{d-m}\kappa$ is an element of $R_{d}$,
whose dimension as a $k$-vector space is $d^{\prime}+1$,
we can write 
$$\theta^{d-m}\kappa 
=\sum_{i=0}^{d^{\prime}} c_{i} x^{s-1-i}y^{t-1-d^{\prime}+i} \>\>\>\>(c_{i} 
\in k).
$$
Then we have $c_{i}+c_{i+1}=0$ for each $i$, 
because	 $\theta \times \theta^{d-m}\kappa =0$ holds.
Hence we find that all $c_{i}$ are non-zero.
Therefore 
$\theta^{d-m} \kappa \times x^{d^{\prime}} = c_{d^{\prime}}x^{s-1}y^{t-1} 
\neq 0$.
Applying Lemma \ref{4}, we finish the proof of this lemma.
\end{proof}

The multiplication map 
$\times \theta^{j} : R_{i} \to R_{i+j}$ is a $k$-linear map.
We denote by $K(i,\,i+j)$ the kernel of this map, and
by $M(i,\,i+j)$ the matrix representation with respect to the canonical bases.

\begin{lem}\label{6}
For each $0 \le i \le s-1$, we have the following:
\begin{enumerate}
\item The map $\times \theta^{t-1-i} : R_{i} \to R_{t-1}$ is injective.
\item The map $\times \theta^{s+t-1-2i} : R_{i} \to R_{s+t-1-i}$ is not injective.
\end{enumerate}
Hence, any non-zero element $\kappa_{i}  \in K(i,\,s+t-1-i) \subseteq R_{i}$ satisfies both
$\theta^{s+t-1-2i}\kappa_{i} = 0$ and $ \theta^{t-1-i} \kappa_{i} \neq 0$.
\end{lem}

\begin{proof}
(1): The map $\times \theta^{t-1-i} : R_{i} \to R_{t-1}$ is represented by the $s \times (i+1)$ matrix:
$$
M(i,\,t-1) =
\begin{pmatrix}
\binom{t-1-i}{t-s} & \binom{t-1-i}{t-s-1} & \cdots & \cdots & \binom{t-1-i}{t-s-i} \\ 
\binom{t-1-i}{t-s+1} & \binom{t-1-i}{t-s} & \cdots & \cdots & \binom{t-1-i}{t-s+1-i} \\ 
\vdots & \vdots &  &  & \vdots \\
\binom{t-1-i}{t-1-i} & \binom{t-1-i}{t-2-i} & \cdots &\cdots  &\binom{t-1-i}{t-1-2i}  \\ 
0&\binom{t-1-i}{t-1-i}  &\cdots  & \cdots &\binom{t-1-i}{t-2i}  \\ 
\vdots& \ddots & \ddots &  &\vdots  \\ 
\vdots&  & \ddots  & \ddots &\vdots  \\ 
0& \cdots & \cdots & 0 & \binom{t-1-i}{t-1-i}
\end{pmatrix}.
$$ 
Hence the map is injective since the rank of $M(i,\,t-1)$ is $i+1$.

(2): It is clear; because 
$i+1=\dim_{k}R_{i} > \dim_{k}R_{s+t-1-i} =i$.
\end{proof}

We now prove Theorem \ref{3}.

\begin{proof}[Proof of Theorem \ref{3}]
We put $n_{0}=n$ and $m_{0}=s+t-2-n_{0}$.
If $m_{0} > 0$, then we have 
$$
\sum_{i_{0}=0}^{m_{0}} k[\theta]x^{i_{0}}=\bigoplus_{i_{0}=0}^{m_{0}} k[\theta]x^{i_{0}} \subseteq R
$$ by Lemma \ref{5}.
If this direct sum coincides with $R$, then we finish the proof.
Suppose not.
By Lemma \ref{6},
we can take an element $\kappa_{(1)} \in K(m_{0}+1,\,n_{0})$,
and then we have $t-1 \le \sigma(k[\theta]\kappa_{(1)}) \le n_{0}-1$.
We put $n_{1} = \sigma(k[\theta]\kappa_{(1)})$ 
and $m_{1}=(n_{0}-1)-n_{1}$.
If $m_{1} > 0$, then we have
$$ 
(\bigoplus_{i_{0}=0}^{m_{0}}k[\theta]x^{i_{0}})+ 
(\sum_{i_{1}=0}^{m_{1}}k[\theta]\kappa_{(1)}x^{i_{1}})
=
\bigoplus_{i_{0}=0}^{m_{0}}k[\theta]x^{i_{0}} \oplus 
\bigoplus_{i_{1}=0}^{m_{1}}k[\theta]\kappa_{(1)}x^{i_{1}} \subseteq R
$$
from Lemma \ref{5}.
Thus, 
we can construct the direct sum of $k[\theta]$-submodules of $R$.
However, since $R$ is finite dimensional, 
this construction will be over in finite steps.
And it is clear that this construction finishes just when $s$-th direct summand is constructed.
By \textit{the Krull-Remak-Schmidt theorem}, 
this decomposition is the indecomposable decomposition of $R$ as a $k[\theta]$-module. 
And this argument does work if some $m_{i}$ is zero.
\begin{figure}
\centering
\includegraphics[clip, scale=0.7]{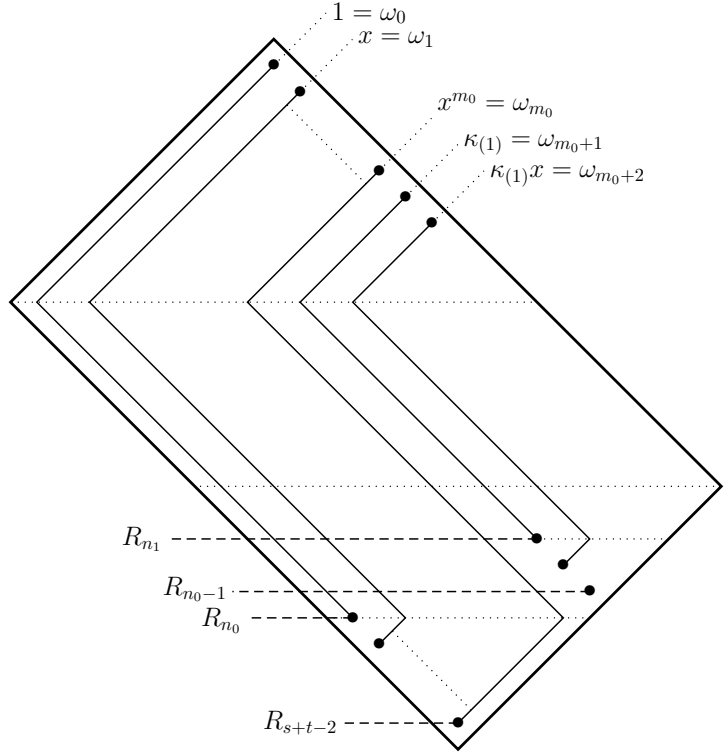}
\caption{Construction of $k[\theta]\omega_{i}$. Each polygonal line stands for $k[\theta]\omega_{i}$.}
\end{figure}
\end{proof}

\begin{rmk}\label{7}
(1)
This proof gives concretely the indecomposable summands of $R$ such as:
\begin{multline*}
k[\theta],\, k[\theta]x,\, \dots,\, k[\theta]x^{m_{0}},\\
k[\theta]\kappa_{(1)},\, k[\theta]\kappa_{(1)}x,\, \dots,\, k[\theta]\kappa_{(1)}x^{m_{1}},\\
\cdots\cdots \\
k[\theta]\kappa_{(r-1)},\, k[\theta]\kappa_{(r-1)}x,\, \dots,\, k[\theta]\kappa_{(r-1)}x^{m_{r-1}},
\end{multline*}
where $\kappa_{(i)}$ means some element in $K(\sum_{j=0}^{i-1}(m_{j}+1),\,n_{i-1})$ and
$m_{i} =(n_{i-1}-1)-n_{i}$, $n_{i} =\sigma(k[\theta]\kappa_{(i)})$. 
Thus, these $\kappa_{(i)},\,m_{i},\,n_{i}$ are determined by the following order:
{\small
\begin{equation*}
(n=)\,n_{0} \to m_{0} \to \kappa_{(1)} \to n_{1} \to m_{1} \to \kappa_{(2)} \to \cdots \to n_{i-1} \to m_{i-1} \to \kappa_{(i)} \to \cdots. 
\end{equation*}}
(Then we define $n_{-1}=s+t-1$, $m_{-1}=0$, and $\kappa_{(0)}=1_{R}$ for convenience).

(2) 
We have to discuss on whether 
the value of $n_{i}=\sigma(k[\theta]\kappa_{(i)})$ depends on 
the choice of an element $\kappa_{(i)} \in K(\sum_{j=0}^{i-1}(m_{j}+1),\,n_{i-1})$.
However, we immediately find that
the numbers $(n_{0},\,n_{1},\dots,\,n_{r-1})$ have to be unique
by the uniqueness of the indecomposable decomposition. 
Therefore we can choose $\kappa_{(i)}$ free. 

(3)
Theorem \ref{3} declares 
the number of Jordan blocks of $J(\alpha,\,s) \otimes J(\beta,\,t)$ is $s$.
\end{rmk}

\begin{dfn}\label{8}
Thus, the particular indecomposable summands 
$$(k[\theta]=)\,k[\theta]\kappa_{(0)}, \, k[\theta]\kappa_{(1)}, \,\dots,\,k[\theta]\kappa_{(r-1)}$$
of $R$ characterize the indecomposable decomposition of $R$. 
So, we shall call each $k[\theta]\kappa_{(i)}$ a \textit{leading module} (of $R$).
And we call 
the number of the indecomposable summands of $R$ 
whose lengths are equal to that of $k[\theta]\kappa_{(i)}$ 
the \textit{leading degree} of $k[\theta]\kappa_{(i)}$.
\end{dfn}

By this result,
if there are $r$ leading modules 
$k[\theta]\kappa_{(0)}, \, k[\theta]\kappa_{(1)}, \,\dots,\,$
$k[\theta]\kappa_{(r-1)}$,
then we have
$$
J(\alpha,\,s)\otimes J(\beta,\, t) = \bigoplus_{i=0}^{r-1}J(\alpha \beta,\, \ell_{i})^{\oplus d_{i}},
$$
where $\ell_{i}$ and $d_{i}$ mean the length and leading degree of $k[\theta]\kappa_{(i)}$ respectively. 

Next, we show a good way to compute a JCF of $J(\alpha,\,s) \otimes J(\beta,\,t)$.
To compute it, we find the lengths and the leading degrees of the leading modules. 

For each $0 \le i \le s-1$, we define a function such as
$$
D_{p}(i) =
\left\{
\begin{array}{cl}
0 & (\text{if the map $\times\theta^{s+t-2-2i} : R_{i} \to R_{s+t-2-i}$ is bijective}) \\
1 & (\text{if the map $\times\theta^{s+t-2-2i} : R_{i} \to R_{s+t-2-i}$ is not bijective}) \\
\end{array} 
\right..
$$
And we put
$$
\varDelta_{p} = (D_{p}(0),\, D_{p}(1),\, \dots,\, D_{p}(s-1)).
$$

\begin{rmk}\label{10}
By Lemma \ref{6} (1),
we have known the map $\times \theta^{t-s} : R_{s-1} \to R_{t-1}$ is always injective (hence, bijective) 
independently of the value of characteristic $p$. 
So $D_{p}(s-1)=0$ holds.
\end{rmk}

By Theorem \ref{3}, we can assume that $R$ is of the form of $\bigoplus_{i=0}^{s-1}k[\theta]\omega_{i}$, i.e. any base of $R$ is that of $\theta^{j}\omega_{i}$. Then the following lemmas hold:

\begin{lem}\label{11}
If an indecomposable summand $k[\theta]\omega_{i}$ is a leading module and $D_{p}(i)=0$. 
Then we have the following:
\begin{enumerate}
\item $\sigma(k[\theta]\omega_{i})=s+t-2-i$. 
Hence the length and the leading degree of $k[\theta]\omega_{i}$ are $s+t-1-2i$ and one respectively.
\item The next indecomposable summand $k[\theta]\omega_{i+1}$ is a leading module.
\end{enumerate}
\end{lem}

\begin{proof}
(1): The map $\times \theta^{s+t-2-2i} : R_{i} \to R_{s+t-2-i}$ is bijective by assumption.
This procedures $\theta^{s+t-2-2i}\omega_{i} \neq 0$ for the generator $\omega_{i}$.
Hence $\sigma(k[\theta]\omega_{i})$ is $s+t-2-i$, and the other statements hold clearly.

(2): It is trivial since the leading degree of $k[\theta]\omega_{i}$ is one.

\begin{figure}
\centering
\includegraphics[clip, scale=0.6]{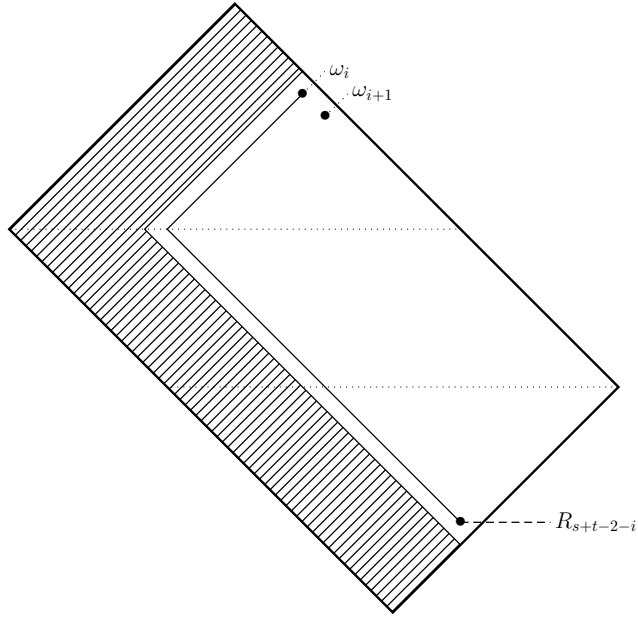}
\caption{The result of Lemma \ref{11}.}
\end{figure}
\end{proof}

\begin{lem}\label{12}
If an indecomposable summand $k[\theta]\omega_{i}$ is a leading module, $D_{p}(i)=D_{p}(i+1)= \cdots = D_{p}(i+f-1)=1\,(f >0)$, and $D_{p}(i+f)=0$. Then we have the following:
\begin{enumerate}
\item $\sigma(k[\theta]\omega_{i})=s+t-2-i-f$. 
Hence the length and the leading degree of $k[\theta]\omega_{i}$ are $s+t-1-2i-f$ and $f+1$ respectively. 
\item The indecomposable summand $k[\theta]\omega_{i+f+1}$ is a leading module.
\end{enumerate}
\end{lem}

\begin{proof}
(1): Put $\nu = \sigma(k[\theta]\omega_{i})$.
Since $D_{p}(i+f)=0$, we have $\theta^{s+t-2-2(i+f)} \times \theta^{f}\omega_{i} = \theta^{s+t-2-2i-f}\omega_{i} \neq 0$.
So $s+t-2-i-f \le \nu \le s+t-2-i$ holds.
Put $\mu=\nu-(s+t-2-i-f)$ and suppose $\mu >0$.
Then
$$
\langle \theta^{s+t-2-(i+f-\mu)}\omega_{0}, \dots, \theta^{s+t-2-2(i+f-\mu)}\omega_{i+f-\mu} \rangle
$$
is a basis of $R_{\nu}$ because the socle of the leading module $K[\theta]\omega_{i}$ is in $R_{\nu}$.
Now $\langle \theta^{i+f-\mu} \omega_{0}, \dots, \omega_{i+f-\mu} \rangle$ is a basis of $R_{i+f-\mu}$.
Hence it is shown that 
\nolinebreak
$\times \theta^{s+t-2-2(i+f-\mu)} : R_{i+f-\mu} \to R_{s+t-2-(i+f-\mu)}$ is bijective.
However, this contradicts $D_{p}(i+f-\mu)=1$.
Therefore $\nu=s+t-2-i-f$. And the other statements hold clearly.

(2): It is trivial from (1).
\end{proof}
 
\begin{figure}
\centering
\includegraphics[clip, scale=0.5]{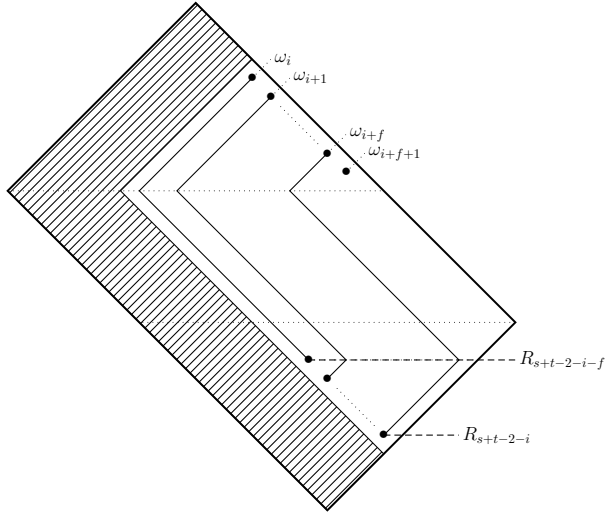}
\caption{The result of Lemma \ref{12}.}
\end{figure}

Since the indecomposable summand $k[\theta]\omega_{0}$ is a leading module,
we can apply Lemma \ref{11} and \ref{12} to the components of $\varDelta_{p}$ inductively.
Thus, via the sequence $\varDelta_{p}$, 
we can compute the lengths and the leading degrees of the leading modules concretely:

\begin{thm}\label{13}
We can compute a JCF of $J(\alpha,\,s) \otimes J(\beta,\,t)$ by using the sequence $\varDelta_{p}$. 
\end{thm}

And we can easily compute the determinant $D(i)$ of the linear map 
$\times \theta^{s+t-2-2i} :R_{i} \to R_{s+t-2-i}$.
In fact, the matrix representation $M(i,\,s+t-2-i)$ is of the form of
$$
\begin{pmatrix}
\binom{s+t-2-2i}{t-1-i}     &\binom{s+t-2-2i}{t-2-i} & \cdots& \cdots &\binom{s+t-2-2i}{t-1-2i} \\
\binom{s+t-2-2i}{t-i} &\binom{s+t-2-2i}{t-1-i}    & \cdots& \cdots &\binom{s+t-2-2i}{t-2i} \\
\vdots&\vdots&\ddots&&\vdots \\
\vdots&\vdots&&\ddots&\vdots \\
\binom{s+t-2-2i}{t-1} &\binom{s+t-2-2i}{t-2} & \cdots& \cdots &\binom{s+t-2-2i}{t-1-i} \\
\end{pmatrix}.
$$
If $p=0$, it is shown by P. C. Roberts \cite{Rob}
that the determinant of the matrix of this form is
computed as follows:
$$
D(i) =\prod_{j=0}^{i} \dfrac{\binom{s+t-2-2i+j}{t-1-i}}{\binom{t-1-i+j}{t-1-i}}. 
$$
And this is true if $p>0$, because $D(i)$ is an integer.

Thus, we get an algorithm for computing a JCF of $J(\alpha,\,s) \otimes J(\beta,\,t)$:

\begin{thm}\label{15}
We are able to compute a JCF of $J(\alpha,\,s) \otimes J(\beta,\,t)$ by taking the following steps:
\begin{enumerate}
\item Computing $D(i)$ for each $0 \le i \le s-1$.
\item Computing the sequence $\varDelta_{p}$. $D_{p}(i)=0$ iff $D(i) \not\equiv 0 \,(\,mod\,\, p\,)$.
\item Applying Theorem \ref{13}.
\end{enumerate} 
\end{thm}

\begin{example}\label{16} 
Let us compute a JCF of $J(\alpha,\,4) \otimes J(\beta,\,5) \,(\alpha \beta \neq 0)$.
The determinants $D(i)$ are
$$
D(0)=\dfrac{\binom{7}{4}}{\binom{4}{4}}=5 \cdot 7, \>\>
D(1)=\dfrac{\binom{5}{3}\binom{6}{3} }{\binom{3}{3}\binom{4}{3}}=2 \cdot 5^{2}, \>\>
D(2)=\dfrac{\binom{3}{2}\binom{4}{2}\binom{5}{2}}{\binom{2}{2}\binom{3}{2}\binom{4}{2}}=2\cdot5, \>\>
D(3)=1.
$$
So the sequence $\varDelta_{p}$ is
$$
\begin{array}{l}
\varDelta_{p}=(0,\,0,\,0,\,0) \>\>(p\neq 2,\,5,\,7),\\
\varDelta_{2}=(0,\,1,\,1,\,0), \\
\varDelta_{5}=(1,\,1,\,1,\,0), \\
\varDelta_{7}=(1,\,0,\,0,\,0). \\
\end{array}
$$
Therefore
{\Small
$$
J(\alpha,\,4) \otimes J(\beta,\, 5) =
\left\{
\begin{array}{ll}
J(\alpha \beta,\, 8) \oplus J(\alpha \beta,\, 6) \oplus J(\alpha \beta,\, 4) \oplus J(\alpha \beta,\, 2) & (p\neq2,\,5,\,7)\\
J(\alpha \beta,\, 8) \oplus J(\alpha \beta,\, 4)^{\oplus 3}  & (p=2)\\
J(\alpha \beta,\,5)^{\oplus 4} & (p=5) \\
J(\alpha \beta,\, 7)^{\oplus 2} \oplus  J(\alpha \beta,\, 4) \oplus J(\alpha \beta,\,2) & (p=7)
\end{array}
\right. .$$
}
\begin{figure}
\centering
\includegraphics[clip, scale=0.75]{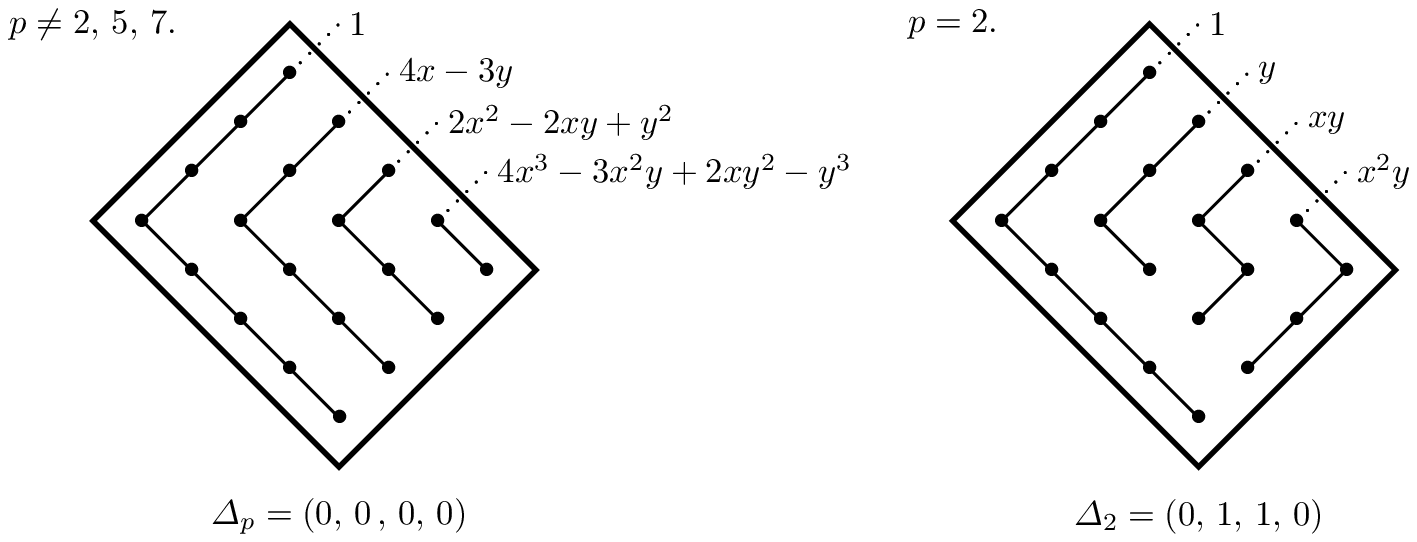}
\centering
\includegraphics[clip, scale=0.75]{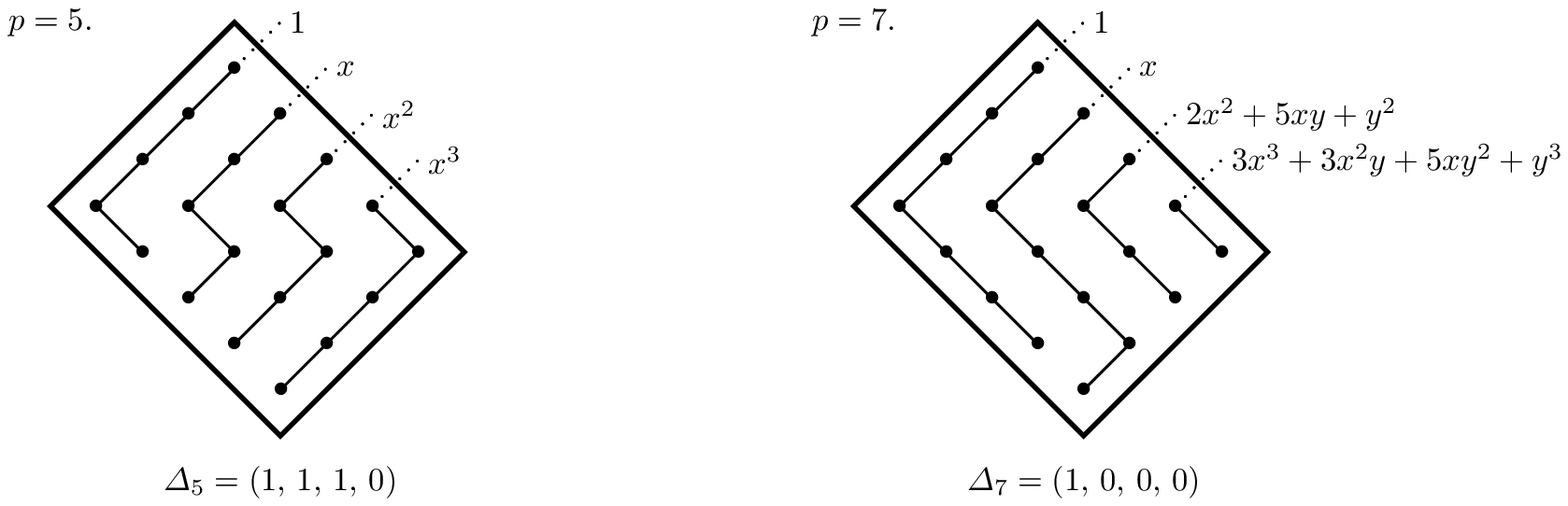}
\caption{The result of Example \ref{16}.}
\end{figure}
\end{example}

If $p=0$ or $p \ge s+t-1$, then the determinants $D(i)$ are clearly all non-zero.
Hence:

\begin{corollary}\label{17} If $p=0$ or $p \ge s+t-1$, then 
$$
J(\alpha,\,s) \otimes J(\beta,\,t) = \bigoplus_{i=0}^{s-1} J(\alpha \beta,\,s+t-1-2i). 
$$
\end{corollary}



{\sc Acknowledgments.}
The authors would like to express their deep gratitude to Yuji Yoshino for suggesting us to explore Jordan decompositions and for a great deal of encouragement as well as making numerous invaluable comments. 
The first author is indebted to Osamu Iyama, Kazuhiko Kurano and Junzo Watanabe for their helpful suggestions. 
The first author would also like to thank Yasuhide Numata and Ryo Takahashi so much for their kind advice. 



\begin{thebibliography}{9999}


\bibitem{BH}
{\sc W.Bruns.}; {\sc J.Herzog.}
Cohen-Macaulay rings. revised edition. Cambridge Studies in Advanced Mathematics, 39. {\it Cambridge University Press, Cambridge}, 1998.

\bibitem{FH}
{\sc W.Fulton.}; {\sc J.Harris.}
Representation theory. A first course. Graduate Texts in Mathematics, 129. Readings in Mathematics. {\it Springer-Verlag, New York}, 1991.

\bibitem{HMNW}
{\sc T.Harima.}; {\sc J.C.Migliore.}; {\sc U.Nagel.}; {\sc J.Watanabe.}
The weak and strong Lefschetz properties for Artinian $K$-algebras.
{\it J. Algebra} {\bf 262} (2003), no. 1, 99--126.

\bibitem{HW}
{\sc T.Harima.}; {\sc J.Watanabe.}
The finite free extension of Artinian $K$-algebras with the strong Lefschetz property.
{\it Rend. Sem. Mat. Univ. Padova} {\bf 110} (2003), 119--146.

\bibitem{H}
{\sc M.Herschend.}
Solution to the Clebsch-Gordan problem for representations of quivers of type $\tilde{\Bbb A}\sb n$.
{\it J. Algebra Appl.} {\bf 4} (2005), no. 5, 481--488.

\bibitem{M}
{\sc I.G.Macdnald.}
Symmetric Functions and Hall Polynomials. second edition. Oxford mathematical monographs. {\it Oxford Science Publications}, 1995.

\bibitem{MV}
{\sc A.Martsinkovsky.}; {\sc A.Vlassov.}
The representation rings of $k[x]$.
preprint, \texttt{http://www.math.neu.edu/martsinkovsky/mathindex.html}.

\bibitem{Rob} 
{\sc P. C. Roberts.} 
A computation of local cohomology.
{\it Contemporary Mathematics.} {\bf 159} (1994), 351--356.

\bibitem{Wak}
{\sc T.Wakamatsu.}
On graded Frobenius algebras.
{\it J. Algebra.} {\bf 267} (2003), 377--395.




\end{thebibliography}
\end{document}